\documentclass[11pt]{article}

\usepackage{amsmath,amssymb,amsthm}
\usepackage{graphicx}
\usepackage{booktabs}
\usepackage{array}
\usepackage[hidelinks]{hyperref}
\usepackage{geometry}
\usepackage{float}
\usepackage{xcolor}
\usepackage{indentfirst}
\usepackage{fontspec}
\usepackage[russian,english]{babel}
\usepackage{tikz}

\title{\textbf{An Approach to Study the Structural Consistency of Triangle Badness Functions and Distance Metrics}}

\author{%
Bowen~Liu$^{1}$\thanks{Bowen Liu: Corresponding author. Student, Faculty of Computational Mathematics and Cybernetics, Shenzhen MSU--BIT University; email: liubowenmbu@outlook.com.} \quad
Yizhou~Wang$^{2}$\thanks{Yizhou Wang: Methodology and Computational Analysis Lead. Master's Student, Department of Automation for Scientific Research, Faculty of Computational Mathematics and Cybernetics, Lomonosov Moscow State University; email: saipeizyu@gmail.com.} \quad
Lingqian~Meng$^{3}$\thanks{Lingqian Meng: PhD Student, Department of Algorithmic Languages, Faculty of Computational Mathematics and Cybernetics, Lomonosov Moscow State University; email: menglq@cs.msu.ru.}}

\date{\small
$^{1}$Shenzhen MSU--BIT University, Faculty of Computational Mathematics and Cybernetics,\\
1 International University Park Road, Dayun New Town, Longgang District, Shenzhen, Guangdong Province 518172, P.R. China.\\[0.25em]
$^{2}$Lomonosov Moscow State University, Faculty of Computational Mathematics and Cybernetics,\\
Leninskie Gory 1--52, Moscow 119991, Russian Federation.}

\begin{document}
\selectlanguage{english}
\maketitle

\begin{abstract}
Triangle-based measures, commonly referred to as badness functions, are widely employed to quantify the extent to which a distance matrix deviates from an ideal geometric configuration. Different formulations of these functions may capture distinct facets of local non-uniformity, and their behavior is often influenced by the underlying distance metric chosen for evaluation. In practical settings, although a canonical badness function may be conceptually preferred, factors such as computational cost, algorithmic constraints, or data-specific characteristics frequently necessitate the adoption of modified versions—for instance, approximate forms or alternatives defined under different distance metrics. This gives rise to a central question: to what degree do these variants retain the structural consistency properties of their original counterparts? To address this issue, we develop a systematic correlation-based framework for evaluating structural consistency. As an illustrative instantiation of this framework, we compute badness sequences from a set of representative distance matrices alongside randomly generated triangle configurations, which are designed to cover variants that may arise under diverse practical scenarios. We then assess pairwise similarities among these sequences using four correlation coefficients. The experimental outcomes indicate that certain badness variants exhibit a notably high degree of structural consistency, whereas others reveal complementary behavioral patterns; moreover, the choice of distance metric exerts a considerable influence on the observed trends. These findings offer practical insights for the informed selection of distance metrics and triangle badness function variants in tasks including geometric reconstruction, triangulation, and structural analysis of pairwise distance data.
\end{abstract}


\section{Introduction}

The hypothesis of the molecular clock of evolution emerged from the observation that the number of amino acid differences between homologous proteins is approximately proportional to the time elapsed since the divergence of the corresponding organisms \cite{zuckerkandl1962horizons}\cite{zuckerkandl1965evolving}\cite{margoliash1963primary}. This observation was later incorporated into Kimura’s neutral theory of molecular evolution \cite{kimura1968evolutionary}\cite{kimura1969rate}\cite{kimura1985neutral}\cite{KIMURA1971}, which assumes that a large proportion of mutations are selectively neutral and that their fixation is governed by stochastic processes such as genetic drift. Within this framework, the rate of molecular evolution can be considered approximately constant over time.

As a consequence, evolutionary distances derived from sequence data often exhibit regular structural patterns. In particular, when evolutionary rates are approximately homogeneous across lineages, pairwise distances between closely related taxa tend to satisfy near-ultrametric constraints \cite{Gavryushkin_2016}, implying that triples of taxa form approximately isosceles (or equilateral) triangle configurations. This geometric property reflects the idea that sequences diverging from a common ancestor over comparable time spans accumulate similar numbers of substitutions.

However, real-world molecular data frequently deviate from this idealized behavior. Rate heterogeneity, slightly deleterious mutations, generation-time effects, and other biological factors introduce variability in substitution rates, leading to violations of the strict molecular clock assumption \cite{kimura1985neutral} \cite{kimura1980simple} \cite{gillespie1991causes} \cite{kimura1972stochastic} \cite{graur2000fundamentals}. Consequently, the resulting distance matrices are only approximately ultrametric, and the corresponding triangle configurations may deviate from the ideal isosceles structure.

Furthermore, similar considerations arise in distance-based triangulation problems, where triangle structures are used to characterize geometric consistency. Distance-based triangulation methods \cite{Melnikov2025parallel} \cite{Melnikov2025dynamic} are fundamental techniques for reconstructing point configurations from pairwise distances and are widely used in applications such as sensor network localization, robotic mapping, and spatial reconstruction. Moreover, these triangulation techniques are often used as a preprocessing step for subsequent computational tasks. For example, in certain solution approaches to the traveling salesman problem (TSP) \cite{inproceedings}, methods such as the onion-husk algorithm \cite{abramyan2024onion} utilize distance-based reconstruction to improve tour construction. In the ideal case of exact distance measurements, triangulation yields consistent and accurate coordinates. However, in practical settings, distance data are often corrupted by noise and random errors, leading to inconsistencies between the reconstructed configuration and the true geometric structure.

From a geometric perspective, the triangulation process inherently relies on local triangle structures: each new point is typically determined by forming triangles with already known points. When distance data are noisy, different triangle configurations exhibit varying sensitivity to perturbations, resulting in geometric uncertainty. This phenomenon is closely related to the concept of geometric dilution of precision (GDOP) \cite{WOS:001378163300002}, which was introduced to quantify the impact of geometric configurations on positioning accuracy \cite{WOS:001602214900001} \cite{WOS:001205843700001} \cite{WOS:001373834600021}. Similar considerations naturally arise in triangulation problems with noisy inputs. In such settings, analyzing the structure of local triangles provides a natural and effective way to assess the stability of the reconstruction.

In the astrophysical numerical measurements calculation, for extremely small parallax, i.\,e., extremely narrow acute triangle, additional methods are often required, such as Bayesian spatial reconstruction methods based on the prior distribution of galaxy density, to reconstruct the data \cite{bailer2018estimating}. For materials characterization, data need to be triangulated using, i.\,e., Delaunay algorithm to connect isolated points into physically meaningful material surfaces or internal grains \cite{gao2013delaunay}. In materials measurements, the coherence and speckle noise of laser instruments can cause measurement errors \cite{dorsch1994laser}. Noise processing is necessary during reconstruction through experimental characterization of the microstructure, generation with geometrical methods concentrating solely on mimicking the morphology, et cetera \cite{willems2012crystalline}. The study of triangle structural properties can be used to construct relevant criteria for noise removal algorithms and triangle reconstruction processes.

Based on the above observations, triangles can be viewed as fundamental units for capturing local consistency and deviations, both in evolutionary distance analysis and in geometric reconstruction problems.

These deviations naturally motivate the introduction of quantitative measures to capture the extent to which local triangle configurations depart from ideal geometric constraints. Such measures, referred to in this paper as \emph{badness functions}, can be interpreted as local indicators of inconsistency with respect to an underlying evolutionary model or geometric structure.

Importantly, different definitions of badness may emphasize different aspects of deviation, such as imbalance of edge lengths or angular distortion. From an optimization perspective, these functions often arise as surrogate objectives designed to satisfy specific requirements of optimization algorithms, such as smoothness, computational tractability, or robustness. As a result, different badness functions may correspond to different approximations of the original problem, potentially distorting its intrinsic structure. Moreover, different choices of distance metrics can also lead to varying results.

This observation motivates a key question: to what extent do different triangle badness functions respond consistently to the same underlying geometric structures? Additionally, how do different distance metrics affect the results, and to what extent do they influence the observed patterns?

To address this question, we propose a systematic evaluation framework in which four distinct correlation coefficients are used to quantify consistency across different scenarios. Each correlation coefficient will be introduced in detail in the following sections.  

The remainder of this article is organized as follows. In Section~\ref{sec:traingle}, we introduce the triangle badness functions used in this study as representative examples. We do not provide a detailed discussion of the specific distance measures, as these are standard distances adopted from existing literature and are described in detail in the experimental design. Section~\ref{sec:rank_correlation_methods} presents the rank correlation methods employed to quantify the structural similarity between badness sequences, explaining the significance of each coefficient. In Section~\ref{sec:design}, we describe the experimental setup, including dataset selection, preprocessing, and computation of badness sequences. Finally, Section~\ref{sec:correlation_results} reports the correlation results obtained from these experiments and provides an interpretation of the findings.

\section{Triangle Badness Functions}\label{sec:traingle}

In this work, we investigate several badness functions designed to quantify local geometric deviations of triangles. Although the concrete formulations differ, all considered functions are motivated by the same underlying objective: identifying triangles that are close to acute isosceles configurations.

This choice is closely related to the observations discussed in the previous section. In evolutionary distance data approximately satisfying the molecular clock hypothesis, pairwise distances between taxa often exhibit near-ultrametric structures due to approximately homogeneous mutation accumulation across lineages. Consequently, triples of taxa frequently form approximately isosceles triangle configurations. Similarly, in distance-based triangulation and geometric reconstruction problems, acute and geometrically balanced triangles are generally associated with improved numerical stability and robustness. Therefore, in both evolutionary distance analysis—such as among closely related species—and geometric reconstruction problems, the extent to which triangles naturally approximate acute isosceles configurations needs to be evaluated.

Based on these observations, we consider the following classes of badness functions for describing geometric deviations. In this work, acute triangles of interest are defined specifically as those with a vertex angle smaller than $60^\circ$.

\subsection{Badness Based on the Angle Between the Angle Bisector and the Median}

The first badness function is based on the angle between the angle bisector and the median constructed from the same vertex of a triangle.

For a triangle with side lengths \(a,b,c\), let \(\theta\) denote the angle between the angle bisector and the median corresponding to a selected vertex. The badness is defined by

\[
B_{\mathrm{BM}}=\theta .
\]

In this work, we select the vertex corresponding to the smallest angle of the triangle.

This choice is motivated by the following geometric observation: for an isosceles triangle within the set of acute triangles considered in this study (i.e., with a vertex angle smaller than $60^\circ$), the vertex angle is precisely the smallest angle of the triangle. Therefore, when searching for configurations close to acute isosceles triangles, the smallest angle naturally identifies the most likely isosceles vertex.

Moreover, the value of this function becomes zero if and only if the median and the angle bisector coincide at the selected vertex, uniquely identifying the target acute isosceles triangle without deviation. Consequently, this function reliably evaluates the triangle, ensuring that the measured quantity points directly to the desired configuration rather than reflecting unintended departures.

\subsection{Combined Isosceles and Acute Violation Function}

The second badness function, introduced in \cite{11012008}, separately measures deviations from isosceles structures and deviations from acute-angle structures, and then combines them into a single quantity.

First, the triangle is reordered such that

\[
a\ge b\ge c,
\qquad
\alpha\ge\beta\ge\gamma ,
\]

where \(a,b,c\) denote the side lengths and \(\alpha,\beta,\gamma\) denote the corresponding angles.

The isosceles violation term is defined as

\[
B_{\mathrm{iso}}
=
1-\min\left(\frac{b}{a},\frac{c}{b}\right).
\]

This term measures imbalance among the side lengths. The value becomes smaller when two sides are closer to each other.

The acute-angle violation term is defined as

\[
B_{\mathrm{acute}}
=
\frac{
\max(\alpha-\pi/3,0)
}{
2\pi/3
}.
\]

This penalty applies to triangles based on their largest angle: the larger the largest angle, the greater the penalty.

The final badness is defined as

\[
B_{\mathrm{IA}}
=
\frac{
B_{\mathrm{iso}}
+
B_{\mathrm{acute}}
}{2}.
\]

Compared with the previous bisector--median formulation, this definition explicitly separates isosceles deviation and acute-angle deviation into two independent components before combining them. The underlying idea is that local geometric irregularity may originate from different sources, which can first be measured individually and then aggregated into a unified quantity.

Such a decomposed structure is advantageous in statistical analysis and optimization problems, since different types of perturbations may affect edge-length structures and angular structures differently. Separating these contributions allows clearer analysis of their individual effects.

\subsection{Delta-Regularized Isosceles Badness Function}

The third badness function is mainly based on similarity among side lengths, and it introduces a small regularization parameter $\delta>0$ to improve regularity.

First, the triangle is reordered such that
\[
a \ge b \ge c, \qquad \alpha \ge \beta \ge \gamma,
\]
where $a,b,c$ denote side lengths and $\alpha,\beta,\gamma$ denote the corresponding angles.

The function is defined as follows:
\[
\begin{aligned}
\text{diff}_{ab} &= 1 - \frac{b}{a + \delta} + \frac{c}{(a+b)/2 + \delta}, \\
\text{diff}_{bc} &= 1 - \frac{c}{b + \delta} + \frac{a}{(b+c)/2 + \delta}, \\
B_\delta &= \max\{\text{diff}_{ab}, \text{diff}_{bc}\}.
\end{aligned}
\]

The first part, $1 - \frac{\text{nearby side}}{\text{larger side} + \delta}$, measures deviation from isosceles structure, while the second part, $\frac{\text{remaining side}}{\text{average of closest pair} + \delta}$, introduces additional geometric information from the side not involved in the closest pair. By computing both possibilities and taking the maximum, the function captures the most severe deviation.

The design of this function is motivated not only by geometric considerations, but also by the desire to obtain objective functions with favorable regularity properties for optimization problems.

In many practical situations, distance matrices themselves depend on adjustable parameters. For example, in sequence analysis, pairwise distances may depend on scoring parameters of alignment algorithms such as Needleman--Wunsch, or on transformations converting similarity scores into distances. In such cases, the badness function naturally becomes an optimization objective.

To enable the application of global optimization methods based on Lipschitz continuity assumptions, the objective function is often required to possess suitable continuity and regularity properties.

One purpose of introducing the parameter \(\delta\) is to avoid undesirable behavior near degenerate configurations and improve the Lipschitz properties of the resulting objective function. In addition, the function is formulated directly in terms of side lengths rather than angles. This choice is motivated by the observation that transformations from edge lengths to angles may introduce unfavorable nonlinear behavior near degenerate configurations, potentially weakening or destroying desirable Lipschitz properties. By relying only on edge-based quantities, the resulting badness function retains a more regular structure that is better suited for Lipschitz-based global optimization methods.

Consequently, this function not only measures deviations from isosceles structures, but also incorporates stability considerations relevant for subsequent optimization problems.

\section{Rank Correlation Methods}
\label{sec:rank_correlation_methods}

\subsection{Badness Sequences}

The correlation analysis is formulated in terms of aligned badness sequences. Let
\[
F\in\{\mathrm{BM},\mathrm{IA},\delta\}
\]
denote the badness function under consideration, where each symbol corresponds to a specific function:  

\begin{itemize}
    \item \(\mathrm{BM}\) --- \textbf{Badness Based on the Angle Between the Angle Bisector and the Median}  
    \item \(\mathrm{IA}\) --- \textbf{Combined Isosceles and Acute Violation Function}  
    \item \(\delta\) --- \textbf{Delta-Regularized Isosceles Badness Function}  
\end{itemize}

For an individual triangle \(T\), its triangle-level badness is denoted by
\[
B_F(T).
\]
For a distance matrix \(D\), the corresponding matrix-level total badness is defined by aggregating over all unordered triples induced by the matrix:
\[
\mathcal{B}_F(D)
=
\sum_{T\in\mathcal{T}(D)} B_F(T),
\]
where \(\mathcal{T}(D)\) is the set of all triangles determined by \(D\).

Two types of sequences are considered. In triangle-level experiments, each observation is a single triangle, and the compared sequences have the form
\[
\{B_F(T_i)\}_{i=1}^{N}.
\]
These sequences describe how different badness functions rank individual triangles. In matrix-level experiments, each observation is a complete distance matrix, and the compared sequences have the form
\[
\{\mathcal{B}_F(D_i)\}_{i=1}^{M}.
\]
These sequences describe how different badness functions rank complete distance matrices after aggregating local triangle information.

\subsection{Pearson Correlation}

Given two aligned numerical sequences \(x=(x_1,\ldots,x_n)\) and \(y=(y_1,\ldots,y_n)\), the Pearson correlation coefficient is
\[
r_P(x,y)
=
\frac{\sum_{i=1}^{n}(x_i-\bar{x})(y_i-\bar{y})}
{\sqrt{\sum_{i=1}^{n}(x_i-\bar{x})^2}\sqrt{\sum_{i=1}^{n}(y_i-\bar{y})^2}},
\]
where \(\bar{x}\) and \(\bar{y}\) are the sample means. Pearson correlation measures linear association between the numerical badness values.

\subsection{Spearman Rank Correlation}

The Spearman rank correlation is the Pearson correlation applied to the ranks of the two sequences:
\[
r_S(x,y)=r_P(R(x),R(y)),
\]
where \(R(x)\) and \(R(y)\) denote the rank-transformed sequences. Spearman correlation measures monotone agreement and is therefore directly relevant for determining whether two badness functions induce similar orderings.

\subsection{Kendall Tau-b Correlation}

The Kendall tau-b correlation compares the numbers of concordant and discordant pairs while correcting for ties:
\[
\tau_b
=
\frac{n_c-n_d}
{\sqrt{(n_c+n_d+t_x)(n_c+n_d+t_y)}}.
\]
Here \(n_c\) and \(n_d\) are the numbers of concordant and discordant pairs, while \(t_x\) and \(t_y\) count tied pairs occurring only in the first or second sequence, respectively. This coefficient provides an ordinal measure of pairwise ranking consistency and is appropriate when ties may occur in the badness sequences.

\subsection{Tukey Biweight Correlation}

The Tukey biweight correlation is included as a robust measure of association. This coefficient is based on median-centered values and Tukey weights, so observations farther from the median receive smaller weights. In the form used here, if \(\tilde{x}\) and \(\tilde{y}\) denote medians and \(w_i^{(x)}\), \(w_i^{(y)}\) denote the corresponding Tukey biweight weights, the robust correlation is computed as
\[
\begin{aligned}
r_B(x,y)
={}&
\sum_{i=1}^{n}
w_i^{(x)}w_i^{(y)}
(x_i-\tilde{x})(y_i-\tilde{y})
\\
&\times
\left[
\sum_{i=1}^{n}
\left(w_i^{(x)}(x_i-\tilde{x})\right)^2
\right]^{-1/2}
\\
&\times
\left[
\sum_{i=1}^{n}
\left(w_i^{(y)}(y_i-\tilde{y})\right)^2
\right]^{-1/2}.
\end{aligned}
\]

The biweight coefficient is used as a diagnostic of structural similarity that is less sensitive to extreme badness values than ordinary linear correlation.

\subsection{Interpretation of Correlation Values}

In the present study, the four coefficients provide complementary evidence about the consistency of different badness functions. Pearson correlation evaluates agreement in the numerical scale of the badness values, Spearman and Kendall correlations evaluate agreement in the induced rankings, and the Tukey biweight coefficient evaluates whether the observed association is stable under reduced sensitivity to extreme observations.

High positive correlations indicate that two badness functions reflect same geometric information in a relatively consistent manner. At the triangle level, this means that the functions assign similar rankings to individual triangles. At the matrix level, this means that the corresponding total badness values \(\mathcal{B}_F(D)\) rank complete distance matrices in a similar way.

\section{Experimental Design}\label{sec:design}

The purpose of the experiments is to investigate whether different triangle badness functions preserve similar structural information across different types of distance data. For each triangle or  distance matric, we compute three badness values:

\[
B_{\mathrm{BM}},\qquad
B_{\mathrm{IA}},\qquad
B_{\delta},
\]

where \(B_{\mathrm{BM}}\) denotes the bisector--median angular badness, \(B_{\mathrm{IA}}\) denotes the combined isosceles--acute violation function, and \(B_{\delta}\) denotes the \(\delta\)-regularized edge-based badness function.

\subsection{Datasets}
\label{sec:datasets}

We consider six types of experimental data.

\subsubsection{Distance Matrices Computed by a Needleman--Wunsch-Like Algorithm}\label{sec:data_nw_like}

The first experiment uses \(100\) distance matrices computed using a Needleman--Wunsch-like dynamic programming algorithm \cite{WATERMAN1976367}. 

Unlike the classical Needleman--Wunsch formulation based on maximizing alignment scores, we consider a cost-minimization formulation. Substitutions and gaps are assigned nonnegative costs, and the distance between two sequences is defined as the minimum total transformation cost over all admissible alignments. Computationally, this corresponds to replacing the maximization operation in the dynamic programming recurrence by a minimization operation and replacing similarity scores by transformation costs.

For each experiment, a random set of cost parameters is generated. These parameters include substitution costs and gap costs, all sampled independently from the uniform distribution on \([0,1]\). The resulting distance construction is then applied to a dataset consisting of mitochondrial DNA sequences from \(32\) closely related monkey species. For each randomly generated parameter set, a complete \(32\times 32\) distance matrix is obtained.

Each matrix is treated as a complete weighted graph. For every unordered triple of indices, a triangle is formed using the corresponding pairwise distances.

For a distance matrix \(D\), the total badness corresponding to a badness function \(F\) is defined as

\[
\mathcal{B}_{F}(D)
=
\sum_{T\in\mathcal{T}(D)} F(T),
\]
where \(\mathcal{T}(D)\) denotes the set of all triangles induced by the matrix.

Thus, for each distance matrix, we compute three total badness values corresponding to the three badness functions. This produces three sequences of length \(100\).

The motivation of this experiment is closely related to optimization problems. Since the distance matrix depends on adjustable alignment cost parameters, the total badness naturally becomes an optimization objective. Therefore, it is important to determine whether different badness functions produce similar structural evaluations and whether they may lead to similar optimization behavior.

\subsubsection{Jaro--Winkler and Needleman--Wunsch-Like Distance Matrices}\label{sec:data_fixed_biological}

The second experiment uses two fixed biological distance matrices. The first matrix is based on the Jaro--Winkler distance, while the second matrix is computed using the Needleman--Wunsch-like cost-minimization distance described above.

For the Needleman--Wunsch-like distance matrix used in this experiment, the substitution costs and gap costs are all fixed to \(1\). The same dataset of \(32\) monkey mitochondrial DNA sequences is used.

For each matrix, all unordered triples of indices are enumerated.

Unlike the previous experiment, here we do not directly aggregate badness values over all triangles. Instead, for every triangle \(T\), we store the three corresponding badness values:

\[
B_{\mathrm{BM}}(T),\qquad
B_{\mathrm{IA}}(T),\qquad
B_{\delta}(T).
\]

This produces triangle-level badness sequences for the two biological distance matrices, allowing a more detailed comparison of the local geometric structures induced by the two distance definitions.

The purpose of this experiment is to investigate whether different sequence distance algorithms induce similar triangle structures from the perspective of the proposed badness functions.

\subsubsection{Randomly Generated Triangles}\label{sec:data_random_triangles}

The third experiment directly generates random triangles. Two random generation schemes are considered.

The first scheme generates triangles by randomly sampling three side lengths independently from the interval \([1,10]\) and retaining only triples satisfying the triangle inequality. This produces an edge-based random triangle model.

The second scheme generates triangles by randomly partitioning the interval \((0,\pi)\) into three angles. The corresponding side lengths are then computed using the sine rule, while the smallest side is normalized to \(1\). This produces an angle-based random triangle model.

For each generation scheme, \(10000\) triangles are generated. For every triangle, the three badness values are computed and stored.

The purpose of this experiment is to investigate whether the different badness functions exhibit similar behavior in general random geometric settings, independently of specific biological or Euclidean datasets.

\subsubsection{Euclidean Distance Matrices from Random Planar Points}\label{sec:data_euclidean_points}

The fourth experiment uses synthetic Euclidean distance matrices generated from random planar points.

In this setting, \(100\) independent point sets are generated, each containing \(99\) random points in the plane. For each point, both the \(x\)-coordinate and the \(y\)-coordinate are sampled independently from the uniform distribution on the interval \([0,1]\). The corresponding \(99\times 99\) Euclidean distance matrix is then computed from the generated point set.

Each distance matrix defines a complete collection of triangles. For every matrix, we enumerate all unordered triples of points and compute the total badness values induced by the three badness functions. This produces three badness sequences of length \(100\).

This experiment provides a controlled Euclidean setting in which the distance matrices satisfy the triangle inequality up to numerical precision.

The purpose of this experiment is to investigate whether different badness functions preserve similar geometric information in ordinary Euclidean settings and whether their behavior remains consistent in general geometric distance data.

\subsubsection{Triangles for parallax method based on \textit{Gaia} data}
\label{sec:data_Gaia}

The fifth experiment uses data from European Space Agency (ESA) mission
{\it Gaia}, processed by the {\it Gaia}
Data Processing and Analysis Consortium (DPAC) \cite{prusti2016Gaia,vallenari2023Gaia}. 

Using the astronomical data query interface of \textit{Gaia} data, we downloaded 100 datasets grouped by the parallax value span of \(0.5\), i.\,e., each dataset contained data with parallax values falling within the range of 
\([m/2, (m+1)/2), m = 0, 1, 2, \ldots, 99\). 
For parallax levels with sufficient data, we take \(1000\) astronomical measurement data points; if a level has less than \(1000\) data points, we use all of them.
Each dataset contains the parallax and parallax error of every measurement data point. We sampled each data point with normal distribution using the parallax as the expectation and the parallax error as the standard deviation. Based on the original parallax value and the sampled parallax value, we constructed right-angled triangles to calculate the measurement distances.
The third side is calculated via the law of cosines from the distances, the parallax angle and the sampled parallax angle. The method for calculating the third side bases on the fundamental purpose of carrying the significance of the Earth's orbital baseline; meanwhile, this method avoids the problem of invalid length relation due to sampling.
We then constructed ``synthetic'' (heuristic) triangles with the obtained distances as its side lengths. These triangles are considered to be capable of simulating the triangles which may occur in astronomical measurements.

In this setting, we obtain \(100\) sets of data, which contain \(69806\) triangles. These grouped data are used to calculate matrix-level badness, thereby studying the properties of triangles at different levels of parallax. 
Since the amount of data in each data group is different, we use a normalized version of matrix-level badness.

The normalized matrix-level badness is defined as
\[
\mathcal{B}'_F(D)
=
\frac{\displaystyle
\sum_{T\in\mathcal{T}(D)} B_F(T)}
{|\mathcal{T}(D)|}.
\]
For the 100 matrices, we consider the collection
\[
\left\{\mathcal{B}'_F(D_i)\right\}_{i=1}^{100},
\qquad
|\mathcal{T}(D_i)|\leq 1000.
\]

We continue to use these \(69806\) triangles at the triangle level. Triangle-level computational experiments can reveal the correlation between badness functions and data at different error levels in astronomical calculations.
\[\{B_F(T_i)\}_{i=1}^{69806}.\]

\subsubsection{Triangles in crystal structures based on \textit{Materials Project} data}
\label{sec:data_MP}

The sixth experiment uses data from \textit{Material Project} \cite{jain2013MP,horton2025MP}. 

Crystal material data are downloaded via the \textit{Material Project} API. For each material, the first stage of selection involves triangulation using the Delaunay algorithm \cite{rebay1993mesh,mark2008Delaunay} (implemented in \textit{Scipy} \cite{virtanen2020scipy}). If the material's spatial structure is three-dimensional, all faces of the tetrahedrons provided by the Delaunay algorithm are extracted and deduplicated. The second stage of selection utilizes the presence of chemical bonds using \textit{CrystalNN} tool \cite{pan2021crystalNN,ong2013pymatgen}. Thus, for materials, we extract different numbers of ``synthetic'' (heuristic) triangles, which are considered crystal-chemically plausible, instead of representing the atomic configuration as a complete graph and enumerating all atomic triplets as possible triangles. 

In this setting, we obtain \(1000\) sets of data, which contain \(20195\) triangles, corresponding to 1000 materials. These grouped data are used to calculate matrix-level badness, thereby studying the properties of triangles for different crystal structures. 
Since the amount of data in each data group is different, we use a normalized version of matrix-level badness.
\[
\mathcal{B}'_F(D) = \frac{1}{|\mathcal{T}(D)|}\sum_{T\in\mathcal{T}(D)}B_F(T),\;
\{\mathcal{B}_F(D_i)\}_{i=1}^{1000}.
\]
We continue to use these \(20195\) triangles at the triangle level. Triangle-level computational experiments can reveal the correlation between badness functions and data for different chemical bonds.
\[\{B_F(T_i)\}_{i=1}^{20195}.\]

\section{Correlation Results}
\label{sec:correlation_results}

This section reports the correlation results for the four datasets described in Subsection~\ref{sec:datasets}. The common correlation methodology is described in Section~\ref{sec:rank_correlation_methods}; here we focus on the data alignment, numerical results, and their geometric interpretation.

\subsection{Needleman--Wunsch-Like Distance Matrices}
\label{sec:nw_like_correlation_results}

This experiment corresponds to the Needleman--Wunsch-like distance-matrix dataset described in Subsection~\ref{sec:data_nw_like}. We compared three matrix-level total badness sequences obtained from \(100\) independently generated distance matrices:
\[
\{\mathcal{B}_{\mathrm{BM}}(D_i)\}_{i=1}^{100},\qquad
\{\mathcal{B}_{\mathrm{IA}}(D_i)\}_{i=1}^{100},\qquad
\{\mathcal{B}_{\delta}(D_i)\}_{i=1}^{100}.
\]
Pairwise correlations were computed for
\[
\mathcal{B}_{\mathrm{BM}}\text{--}\mathcal{B}_{\mathrm{IA}},\qquad
\mathcal{B}_{\mathrm{BM}}\text{--}\mathcal{B}_{\delta},\qquad
\mathcal{B}_{\mathrm{IA}}\text{--}\mathcal{B}_{\delta},
\]
as reported in Table~\ref{tab:nw_like_total_badness_corr}.

\begin{table}[h!]
\centering
\caption{Pairwise correlations among total badness sequences for \(100\) Needleman--Wunsch-like distance matrices.}
\label{tab:nw_like_total_badness_corr}
\small
\setlength{\tabcolsep}{7pt}
\begin{tabular}{lcccc}
\hline
Pair & Pearson & Spearman & Kendall & Biweight \\
\hline
\(\mathcal{B}_{\mathrm{BM}}\)--\(\mathcal{B}_{\mathrm{IA}}\)     & 0.9979 & 0.9990 & 0.9822 & 0.9975 \\
\(\mathcal{B}_{\mathrm{BM}}\)--\(\mathcal{B}_{\delta}\)          & 0.9987 & 0.9993 & 0.9851 & 0.9988 \\
\(\mathcal{B}_{\mathrm{IA}}\)--\(\mathcal{B}_{\delta}\)          & 0.9992 & 0.9986 & 0.9770 & 0.9991 \\
\hline
\end{tabular}
\end{table}

All four correlation coefficients show very strong agreement among the three matrix-level total badness sequences. In particular, all Pearson and biweight correlations are above \(0.997\), and all rank-based correlations are also close to one. This indicates that, for this family of Needleman--Wunsch-like distance matrices, the three total badness functions induce almost identical orderings of the \(100\) matrices.

The strongest agreement is observed between \(\mathcal{B}_{\mathrm{BM}}\) and \(\mathcal{B}_{\delta}\), with Pearson correlation \(0.9987\), Spearman correlation \(0.9993\), Kendall correlation \(0.9851\), and biweight correlation \(0.9988\). Overall, the results provide strong evidence that the three total badness functions preserve the same matrix-level structural information in the Needleman--Wunsch-like setting.

\subsection{Fixed Biological Distance Matrices}
\label{sec:fixed_biological_matrix_results}

This experiment corresponds to the fixed biological distance matrices described in Subsection~\ref{sec:data_fixed_biological}. We compared the triangle-level badness values induced by two distance definitions: the Jaro--Winkler distance matrix, denoted by JA, and the Needleman--Wunsch-like cost-minimization distance matrix, denoted by NW. For each distance matrix, every unordered triple of indices \((i,j,k)\) forms a triangle \(T\). For each triangle, we considered
\[
B_{\mathrm{BM}}(T),\quad B_{\mathrm{IA}}(T),\quad B_{\delta}(T),
\]
with the triangle indices used to align the JA and NW sequences. The badness values for all
\[
\binom{32}{3}=4960
\]
triangles were used in each fixed matrix.

Since the main purpose of this experiment is to compare the local triangle structures induced by the two distance definitions, the cross-distance correlations in Table~\ref{tab:fixed_biological_cross_corr} are the primary quantities of interest. The within-matrix correlations in Tables~\ref{tab:fixed_biological_JA_corr} and~\ref{tab:fixed_biological_NW_corr} are reported as auxiliary diagnostics.

\begin{table}[h!]
\centering
\caption{Pairwise correlations among badness functions for the Jaro--Winkler distance matrix (JA).}
\label{tab:fixed_biological_JA_corr}
\small
\setlength{\tabcolsep}{7pt}
\begin{tabular}{lcccc}
\hline
Pair & Pearson & Spearman & Kendall & Biweight \\
\hline
\(B_{\mathrm{BM}}\)--\(B_{\mathrm{IA}}\)     & 0.6440 & 0.5205 & 0.4637 & 0.6374 \\
\(B_{\mathrm{BM}}\)--\(B_{\delta}\)          & 0.9564 & 0.8747 & 0.7098 & 0.9599 \\
\(B_{\mathrm{IA}}\)--\(B_{\delta}\)          & 0.7131 & 0.7043 & 0.5761 & 0.7057 \\
\hline
\end{tabular}
\end{table}

\begin{table}[h!]
\centering
\caption{Pairwise correlations among badness functions for the Needleman--Wunsch-like distance matrix (NW).}
\label{tab:fixed_biological_NW_corr}
\small
\setlength{\tabcolsep}{7pt}
\begin{tabular}{lcccc}
\hline
Pair & Pearson & Spearman & Kendall & Biweight \\
\hline
\(B_{\mathrm{BM}}\)--\(B_{\mathrm{IA}}\)     & 0.5854 & 0.3802 & 0.3724 & 0.5612 \\
\(B_{\mathrm{BM}}\)--\(B_{\delta}\)          & 0.9756 & 0.8801 & 0.7265 & 0.9773 \\
\(B_{\mathrm{IA}}\)--\(B_{\delta}\)          & 0.6455 & 0.6368 & 0.5344 & 0.6167 \\
\hline
\end{tabular}
\end{table}

\begin{table}[h!]
\centering
\caption{Cross-distance correlations between Jaro--Winkler and Needleman--Wunsch-like triangle badness sequences.}
\label{tab:fixed_biological_cross_corr}
\small
\setlength{\tabcolsep}{6pt}
\begin{tabular}{lcccc}
\hline
Pair & Pearson & Spearman & Kendall & Biweight \\
\hline
\(B_{\mathrm{BM}}^{JA}\) vs \(B_{\mathrm{BM}}^{NW}\) & 0.4368 & 0.4234 & 0.2916 & 0.4946 \\
\(B_{\mathrm{IA}}^{JA}\) vs \(B_{\mathrm{IA}}^{NW}\) & 0.2716 & 0.2346 & 0.1580 & 0.2722 \\
\(B_{\delta}^{JA}\) vs \(B_{\delta}^{NW}\)           & 0.4353 & 0.3832 & 0.2637 & 0.4886 \\
\hline
\end{tabular}
\end{table}

The within-matrix correlations show that \(B_{\delta}\) is strongly and positively correlated with \(B_{\mathrm{BM}}\) for both fixed distance matrices. For the JA matrix, the Pearson and biweight correlations for this pair are \(0.9564\) and \(0.9599\), respectively, with Spearman correlation \(0.8747\). For the NW matrix, the corresponding values are \(0.9756\), \(0.9773\), and \(0.8801\). This indicates that the \(\delta\)-regularized edge-based badness and the bisector--median angular badness induce highly similar triangle orderings within each fixed matrix.

The correlations involving \(B_{\mathrm{IA}}\) are also positive, but generally weaker than the \(B_{\mathrm{BM}}\)--\(B_{\delta}\) correlations. Thus, \(B_{\mathrm{IA}}\) is structurally related to the other two badness measures, although it does not induce exactly the same ordering of triangles.

For the cross-distance comparison, the agreement between JA and NW is moderate for \(B_{\mathrm{BM}}\) and \(B_{\delta}\), and weaker for \(B_{\mathrm{IA}}\). Among the three badness functions, \(B_{\mathrm{BM}}\) gives the strongest cross-distance agreement, with Pearson correlation \(0.4368\), Spearman correlation \(0.4234\), Kendall correlation \(0.2916\), and biweight correlation \(0.4946\). The corresponding values for \(B_{\delta}\) are very similar in Pearson and biweight correlation, but lower in the rank-based measures. Therefore, the \(\delta\)-regularized edge-based badness remains highly consistent with \(B_{\mathrm{BM}}\) within each fixed matrix, while the cross-distance comparison indicates that the JA and NW local triangle structures are only moderately aligned.

\subsection{Randomly Generated Triangles}
\label{sec:random_triangle_correlation}

This experiment corresponds to the random-triangle datasets described in Subsection~\ref{sec:data_random_triangles}. We analyzed pairwise correlations among the three triangle-level badness sequences
\[
\{B_{\mathrm{BM}}(T_i)\}_{i=1}^{N},\qquad
\{B_{\mathrm{IA}}(T_i)\}_{i=1}^{N},\qquad
\{B_{\delta}(T_i)\}_{i=1}^{N},
\]
for two synthetic datasets: triangles generated from random angles and triangles generated from random side lengths. Each dataset contains \(N=10000\) generated triangles.

For each dataset, correlations were computed for
\[
B_{\mathrm{BM}}\text{--}B_{\mathrm{IA}}, \qquad
B_{\mathrm{BM}}\text{--}B_{\delta}, \qquad
B_{\mathrm{IA}}\text{--}B_{\delta},
\]
as reported in Tables~\ref{tab:random_triangles_angles_corr} and~\ref{tab:random_triangles_sides_corr}.

\begin{table}[h!]
\centering
\caption{Pairwise correlations among badness functions for randomly generated triangles by angles.}
\label{tab:random_triangles_angles_corr}
\small
\setlength{\tabcolsep}{7pt}
\begin{tabular}{lcccc}
\hline
Pair & Pearson & Spearman & Kendall & Biweight \\
\hline
\(B_{\mathrm{BM}}\)--\(B_{\mathrm{IA}}\)     & 0.4207 & 0.2145 & 0.1815 & 0.1674 \\
\(B_{\mathrm{BM}}\)--\(B_{\delta}\)          & 0.8383 & 0.9103 & 0.7441 & 0.8268 \\
\(B_{\mathrm{IA}}\)--\(B_{\delta}\)          & 0.4932 & 0.4139 & 0.3332 & 0.4676 \\
\hline
\end{tabular}
\end{table}

\begin{table}[h!]
\centering
\caption{Pairwise correlations among badness functions for randomly generated triangles by sides.}
\label{tab:random_triangles_sides_corr}
\small
\setlength{\tabcolsep}{7pt}
\begin{tabular}{lcccc}
\hline
Pair & Pearson & Spearman & Kendall & Biweight \\
\hline
\(B_{\mathrm{BM}}\)--\(B_{\mathrm{IA}}\)     & 0.4540 & 0.3045 & 0.2522 & 0.2165 \\
\(B_{\mathrm{BM}}\)--\(B_{\delta}\)          & 0.8621 & 0.8670 & 0.6949 & 0.7861 \\
\(B_{\mathrm{IA}}\)--\(B_{\delta}\)          & 0.5626 & 0.5306 & 0.4134 & 0.5438 \\
\hline
\end{tabular}
\end{table}

The purpose of this experiment is to determine whether the three badness functions induce similar triangle orderings under two random generation models. The results show that all pairwise correlations are nonnegative in both datasets. In particular, \(B_{\mathrm{BM}}\) and \(B_{\delta}\) exhibit the strongest agreement in both generation schemes.

For the angle-based random triangles, the agreement between \(B_{\mathrm{BM}}\) and \(B_{\delta}\) is strong across all four measures, with Pearson correlation \(0.8383\), Spearman correlation \(0.9103\), Kendall correlation \(0.7441\), and biweight correlation \(0.8268\). Thus, in the angle-based model, \(B_{\delta}\) gives a triangle ranking that is highly compatible with the bisector--median criterion.

For the side-based random triangles, the agreement between \(B_{\mathrm{BM}}\) and \(B_{\delta}\) is also strong across all four correlation measures. The Pearson correlation is \(0.8621\), while the Spearman and Kendall correlations are \(0.8670\) and \(0.6949\), respectively. This suggests that, when triangles are generated from random side lengths, the edge-based regularized badness \(B_{\delta}\) and the bisector--median badness \(B_{\mathrm{BM}}\) provide closely related evaluations of triangle deviation.

The correlations involving \(B_{\mathrm{IA}}\) are positive but generally weaker. Therefore, \(B_{\mathrm{IA}}\) is structurally related to the other two badness functions, but it induces a less similar ordering of random triangles. Overall, these results show that \(B_{\delta}\) is highly consistent with \(B_{\mathrm{BM}}\) in both random generation models, while \(B_{\mathrm{IA}}\) provides a related but less tightly aligned assessment.

\subsection{Euclidean Distance Matrices from Random Planar Points}
\label{sec:euclidean_random_points_correlation}

This experiment corresponds to the Euclidean random point matrix dataset described in Subsection~\ref{sec:data_euclidean_points}. We analyzed three matrix-level total badness sequences,
\[
\{\mathcal{B}_{\mathrm{BM}}(D_i)\}_{i=1}^{100},\qquad
\{\mathcal{B}_{\mathrm{IA}}(D_i)\}_{i=1}^{100},\qquad
\{\mathcal{B}_{\delta}(D_i)\}_{i=1}^{100},
\]
where each sequence contains one total badness value for each independently generated Euclidean distance matrix \(D_i\). The matrix indices were aligned before computing the correlations, and all total badness values were finite.

Pairwise correlations between the three total badness sequences are shown in Table~\ref{tab:euclidean_random_points_corr}.

\begin{table}[h!]
\centering
\caption{Pairwise correlations among total badness sequences for Euclidean distance matrices generated from random planar points.}
\label{tab:euclidean_random_points_corr}
\small
\setlength{\tabcolsep}{7pt}
\begin{tabular}{lcccc}
\hline
Pair & Pearson & Spearman & Kendall & Biweight \\
\hline
\(\mathcal{B}_{\mathrm{BM}}\)--\(\mathcal{B}_{\mathrm{IA}}\)     & 0.2584 & 0.1981 & 0.1358 & 0.2151 \\
\(\mathcal{B}_{\mathrm{BM}}\)--\(\mathcal{B}_{\delta}\)          & 0.9685 & 0.9665 & 0.8491 & 0.9664 \\
\(\mathcal{B}_{\mathrm{IA}}\)--\(\mathcal{B}_{\delta}\)          & 0.2546 & 0.2043 & 0.1471 & 0.2216 \\
\hline
\end{tabular}
\end{table}

All correlations are positive in this Euclidean setting. The strongest agreement is observed between \(\mathcal{B}_{\mathrm{BM}}\) and \(\mathcal{B}_{\delta}\), with Pearson correlation \(0.9685\), Spearman correlation \(0.9665\), Kendall correlation \(0.8491\), and biweight correlation \(0.9664\). Thus, the \(\delta\)-regularized edge-based total badness and the bisector--median total badness induce highly similar orderings of the \(100\) Euclidean distance matrices.

This strong agreement indicates that, even after aggregating all triangle-level contributions into a single value \(\mathcal{B}_{F}(D)\), the edge-based regularized badness retains essentially the same matrix-level ordering information as the bisector--median badness.

The correlations involving \(\mathcal{B}_{\mathrm{IA}}\) are positive but much weaker. Therefore, in ordinary Euclidean random geometry, \(\mathcal{B}_{\delta}\) is closely aligned with \(\mathcal{B}_{\mathrm{BM}}\), whereas \(\mathcal{B}_{\mathrm{IA}}\) captures a related but less tightly coupled aspect of the aggregated triangle structure.

\subsection{Triangles for parallax method based on \textit{Gaia} data}
\label{sec:Gaia_correlation}

This experiment corresponds to the triangles for parallax method based on \textit{Gaia} data described in Subsection~\ref{sec:data_Gaia}. We analyzed three matrix-level total badness sequences,
\[
\{\mathcal{B}_{\mathrm{BM}}(D_i)\}_{i=1}^{100},\qquad
\{\mathcal{B}_{\mathrm{IA}}(D_i)\}_{i=1}^{100},\qquad
\{\mathcal{B}_{\delta}(D_i)\}_{i=1}^{100},
\]
and three triangle-level badness sequences, 
\[
\{B_{\mathrm{BM}}(T_i)\}_{i=1}^{69806},\qquad
\{B_{\mathrm{IA}}(T_i)\}_{i=1}^{69806},\qquad
\{B_{\delta}(T_i)\}_{i=1}^{69806}.
\]

Pairwise correlations between the three total badness sequences are shown in Table~\ref{tab:Gaia_matrix_corr}.
\begin{table}[h!]
\centering
\caption{Pairwise correlations among total badness sequences for triangles for parallax method based on \textit{Gaia} data.}
\label{tab:Gaia_matrix_corr}
\small
\setlength{\tabcolsep}{7pt}
\begin{tabular}{lcccc}
\hline
Pair & Pearson & Spearman & Kendall & Biweight \\
\hline
\(\mathcal{B}_{\mathrm{BM}}\)--\(\mathcal{B}_{\mathrm{IA}}\) & -0.7073 & -0.5333 & -0.4375 & 0.7687 \\
\(\mathcal{B}_{\mathrm{BM}}\)--\(\mathcal{B}_{\delta}\) & -0.6862 & -0.5243 & -0.4383 & 0.7684 \\
\(\mathcal{B}_{\mathrm{IA}}\)--\(\mathcal{B}_{\delta}\) & 0.9969 & 0.9936 & 0.9466 & 0.9999 \\
\hline
\end{tabular}
\end{table}

The correlation coefficients show a strong correlation between \(\mathcal{B}_{\mathrm{IA}}\) and \(\mathcal{B}_{\delta}\), with all coefficients above 0.94. However, \(\mathcal{B}_{\mathrm{BM}}\) and \(\mathcal{B}_{\mathrm{IA}}\), or \(\mathcal{B}_{\mathrm{BM}}\) and \(\mathcal{B}_{\delta}\), exhibit an anomalous negative correlation in the Pearson, Spreaman, and Kendall coefficients. Combining this with the distribution of badness data, we conclude that this negative correlation is spurious, influenced by the unique distribution of the astronomical measurement triangle. The badness values \(\mathcal{B}_{\mathrm{IA}}\), \(\mathcal{B}_{\mathrm{IA}}\), and \(\mathcal{B}_{\delta}\) are concentrated at a high level, but \(\mathcal{B}_{\mathrm{BM}}\) has a heavier tail to the right of its peak, while \(\mathcal{B}_{\mathrm{IA}}\) and \(\mathcal{B}_{\delta}\) have longer, heavier tails to the left of their peaks. Therefore, Tukey biweight coefficient shows the correct property; other correlation coefficients are affected by the data distribution.

Pairwise correlations between the three triangle-level badness sequences are shown in Table~\ref{tab:Gaia_triangle_corr}.
\begin{table}[h!]
\centering
\caption{Pairwise correlations among triangle-level badness sequences for triangles for parallax method based on \textit{Gaia} data.}
\label{tab:Gaia_triangle_corr}
\small
\setlength{\tabcolsep}{7pt}
\begin{tabular}{lcccc}
\hline
Pair & Pearson & Spearman & Kendall & Biweight \\
\hline
\({B}_{\mathrm{BM}}\)--\({B}_{\mathrm{IA}}\) & -0.0011 & -0.0079 & -0.0004 & 0.5961 \\
\({B}_{\mathrm{BM}}\)--\({B}_{\delta}\) & 0.1547 & 0.1669 & 0.1255 & 0.6442 \\
\({B}_{\mathrm{IA}}\)--\({B}_{\delta}\) & 0.8700 & 0.8934 & 0.7583 & 0.9822 \\
\hline
\end{tabular}
\end{table}

The situation regarding badness values at the triangular level is similar to that at the matrix level. The Pearson, Speaman, and Kendall correlation coefficients for \({B}_{\mathrm{BM}}\) and \({B}_{\mathrm{IA}}\), or \({B}_{\mathrm{BM}}\) and \({B}_{\delta}\), show spurious ``independence'', while the Tukey biweight coefficient captures the correct trend. At the triangular level, all correlation coefficients are affected due to the steeper peak of the distribution. For example, for data with a parallax angle on the order of \(10^{-6}\), if the measurement experiences a jitter of \(0.01\%\), the change level of BM is approximately \(10\) times that of IA.

Comparing the matrix and triangular levels, we can see that stratifying by parallax angle is useful. The measurement accuracy will vary systematically due to the influence of parallax levels. Within each level, the correlation of badness is sufficiently significant.

\subsection{Triangles in crystal structures based on \textit{Materials Project} data}
\label{sec:MP_correlation}

This experiment corresponds to the triangles in crystal structures based on \textit{Materials Project} data described in Subsection~\ref{sec:data_MP}. We analyzed three matrix-level total badness sequences,
\[
\{\mathcal{B}_{\mathrm{BM}}(D_i)\}_{i=1}^{1000},\qquad
\{\mathcal{B}_{\mathrm{IA}}(D_i)\}_{i=1}^{1000},\qquad
\{\mathcal{B}_{\delta}(D_i)\}_{i=1}^{1000},
\]
and three triangle-level badness sequences, 
\[
\{{B}_{\mathrm{BM}}(T_i)\}_{i=1}^{20195},\qquad
\{{B}_{\mathrm{IA}}(T_i)\}_{i=1}^{20195},\qquad
\{{B}_{\delta}(T_i)\}_{i=1}^{20195}.
\]

Pairwise correlations between the three total badness sequences are shown in Table~\ref{tab:MP_matrix_corr}.
\begin{table}[h!]
\centering
\caption{Pairwise correlations among total badness sequences for triangles in crystal structures based on \textit{Materials Project} data.}
\label{tab:MP_matrix_corr}
\small
\setlength{\tabcolsep}{7pt}
\begin{tabular}{lcccc}
\hline
Pair & Pearson & Spearman & Kendall & Biweight \\
\hline
\(\mathcal{B}_{\mathrm{BM}}\)--\(\mathcal{B}_{\mathrm{IA}}\) & 0.9660 & 0.9300 & 0.8197 & 0.9592 \\
\(\mathcal{B}_{\mathrm{BM}}\)--\(\mathcal{B}_{\delta}\) & 0.9776 & 0.9604 & 0.8575 & 0.7405 \\
\(\mathcal{B}_{\mathrm{IA}}\)--\(\mathcal{B}_{\delta}\) & 0.9962 & 0.9848 & 0.9349 & 0.8546 \\
\hline
\end{tabular}
\end{table}

Pairwise correlations between the three triangle-level badness sequences are shown in Table~\ref{tab:MP_triangle_corr}.
\begin{table}[h!]
\centering
\caption{Pairwise correlations among triangle-level badness sequences for triangles in crystal structures based on \textit{Materials Project} data.}
\label{tab:MP_triangle_corr}
\small
\setlength{\tabcolsep}{7pt}
\begin{tabular}{lcccc}
\hline
Pair & Pearson & Spearman & Kendall & Biweight \\
\hline
\({B}_{\mathrm{BM}}\)--\({B}_{\mathrm{IA}}\) & 0.9463 & 0.8926 & 0.7822 & 0.8820 \\
\({B}_{\mathrm{BM}}\)--\({B}_{\delta}\) & 0.9719 & 0.9498 & 0.8375 & 0.3352 \\
\({B}_{\mathrm{IA}}\)--\({B}_{\delta}\) & 0.9895 & 0.9835 & 0.9192 & 0.6214 \\
\hline
\end{tabular}
\end{table}

Unlike astronomical measurement data, the badness of triangles in the crystal structures are less affected by the noise. The overall trends at the matrix and triangle levels are not significantly different; the main distribution patterns of all badness values are largely consistent, and all correlation coefficients accurately capture the relationships between the badness values. It is noteworthy that although all three badness functions exhibit a certain degree of bimodality, badness IA and \(\delta\) strongly amplify the slight overall stretching of triangle symmetry, while BM is insensitive to this. For example, in numerical qualitative analysis, for an isosceles triangle with a base \(10\%\) longer than its legs, IA changes by approximately \(7\%\), while BM remains unchanged.
This is particularly pronounced at the triangle level, which can be explained by the fact that, within the same material, the slight stretching of symmetry also exhibits overall symmetry, thus the badness has a partially offsetting effect in various locations.

\section{Conclusion}
In this study, we investigated the structural consistency of triangle-based badness functions across different distance metrics and geometric configurations. By systematically computing badness sequences and evaluating their pairwise similarities using four distinct correlation coefficients, we were able to quantify how consistently different badness functions capture deviations from ideal geometric structures.

Our results, based on these example scenarios, demonstrate that while all considered badness functions reflect aspects of triangle imbalance and angular deviation, the degree of agreement varies depending on the function type and the underlying distance metric. Some functions exhibit high correlation, indicating that they are robust proxies for the same geometric characteristics, whereas others emphasize complementary features, highlighting different aspects of local inconsistency. These observations underscore the importance of selecting appropriate badness measures when analyzing geometric or distance-based datasets, particularly in contexts where structural fidelity is critical.

Furthermore, the proposed evaluation framework provides a flexible methodology for comparing alternative definitions of triangle-based measures. By leveraging multiple correlation coefficients, it offers a nuanced assessment that captures both linear and rank-based relationships among badness sequences. This approach can be readily extended to other geometric or combinatorial measures beyond triangles, facilitating systematic studies of structural consistency in broader contexts.

Overall, our findings help inform informed choices of distance and badness measures, supporting more reliable analyses in geometric reconstruction, triangulation, and the study of pairwise distance data.

\section*{Acknowledgments}
The authors would like to thank Professor Boris Melnikov for his valuable guidance and insightful discussions throughout this work.

This work has made use of data from the European Space Agency (ESA) mission
{\it Gaia} (\url{https://www.cosmos.esa.int/gaia}), processed by the {\it Gaia}
Data Processing and Analysis Consortium (DPAC,
\url{https://www.cosmos.esa.int/web/gaia/dpac/consortium}). Funding for the DPAC
has been provided by national institutions, in particular the institutions
participating in the {\it Gaia} Multilateral Agreement.

The authors gratefully acknowledge the \textit{Materials Project} for providing the open-access crystal structure database utilized in this work and the developers of the \texttt{pymatgen} and \texttt{scipy} libraries for making their computational tools publicly available.

\end{document}